\documentclass[12pt]{amsart}

\usepackage{amsmath,amssymb,amscd,amsthm,epsf,graphicx}

\numberwithin{equation}{section}
\theoremstyle{plain}
\newtheorem{theorem}[equation]{Theorem}

\newtheorem{lemma}[equation]{Lemma}
\newtheorem{corollary}[equation]{Corollary}

\theoremstyle{remark}
\newtheorem{remark}[equation]{Remark}

\theoremstyle{definition}
\newtheorem{definition}[equation]{Definition}

\newcommand{\lra}{\longrightarrow}
\newcommand{\ra}{\rightarrow}
\newcommand{\restr}{\mbox{\Large \(|\)\normalsize}}

\newcommand{\B}{{\mathcal B}}

\newcommand{\D}{{\mathcal D}}

\newcommand{\G}{{\mathcal G}}

\newcommand{\h}{{\mathcal H}}

\newcommand{\N}{\mathbb N}

\newcommand{\R}{\mathbb R}

\newcommand{\U}{{\mathcal U}}
\newcommand{\V}{{\mathcal V}}

\newcommand{\Z}{\mathbb Z}

\newcommand{\acts}{\curvearrowright}

\newcommand{\isom}{\operatorname{Isom}}

\renewcommand{\span}{\operatorname{span}}

\def\D{\partial}
\newcommand{\al}{\alpha}
\def\de{\delta}

\def\ga{\gamma}
\def\Ga{\Gamma}

\def\la{\lambda}

\def\lra{\longrightarrow}

\def\om{\omega}

\def\ra{\rightarrow}

\def\be{\beta}

\def\defeq{:=}

\def\XXint#1#2#3{{\setbox0=\hbox{$#1{#2#3}{\int}$}
     \vcenter{\hbox{$#2#3$}}\kern-.5\wd0}}

\begin{document}

\title[Groups of polynomial growth]{A new proof of Gromov's theorem 
on groups of polynomial growth}
\author{Bruce Kleiner}
\thanks{Supported by NSF Grant DMS 0701515}
\date{\today}
\maketitle

\begin{abstract}
We give a proof of Gromov's theorem that any finitely 
generated group of polynomial growth has a finite index 
nilpotent subgroup.  The proof does not rely on the 
Montgomery-Zippin-Yamabe structure theory of locally 
compact groups.
\end{abstract}

\tableofcontents

\section{Introduction}

\subsection{Statement of results}
\begin{definition}
Let $G$ be a finitely generated group,
and let $B_G(r)\subset G$ denote the ball centered at 
$e\in G$ with respect to some fixed word norm on $G$. 
The group $G$  has {\bf polynomial growth} if 
for some $d\in (0,\infty)$
\begin{equation}
\label{eqnstrongpolygrowth}
\limsup_{r\ra\infty}\;\frac{|B_G(r)|}{r^d}<\infty,
\end{equation}
and  has {\bf weakly polynomial growth} if 
for some $d\in (0,\infty)$
\begin{equation}
\label{eqnweakpolygrowth}
\liminf_{r\ra\infty}\;\frac{|B_G(r)|}{r^d}<\infty,
\end{equation}
\end{definition}

We give a  proof of the following  special case
of a theorem of Colding-Minicozzi, without using
Gromov's theorem on groups of polynomial growth:

\begin{theorem}[\cite{coldingminicozzi}]
\label{thmfinitedim}
Let $\Ga$ be a Cayley graph of a group $G$ of weakly polynomial 
growth,
and $d\in [0,\infty)$.
Then the space of harmonic functions on $\Ga$ with polynomial growth
at most $d$ is finite dimensional.
\end{theorem}
Note that although \cite{coldingminicozzi} stated the result for
groups of polynomial growth, their proof also works 
for groups of weakly polynomial growth, in view of \cite{wilkie}.

We then use this to derive the following corollaries:
\begin{corollary}
\label{corinfiniterep}
If $G$ is an infinite group of weakly polynomial growth, then 
$G$ admits a finite dimensional linear representation $G\ra GL(n,\R)$
with  infinite image.
\end{corollary}

\begin{corollary}[ \cite{polygrowth,wilkie}]
\label{corgromov}
If $G$ is a  group with weakly polynomial growth, then $G$
is virtually nilpotent.
\end{corollary}

We emphasize that our proof
of Corollary \ref{corgromov} yields a new proof of Gromov's theorem on 
groups of polynomial growth,  which does not involve the 
Montgomery-Zippin-Yamabe structure theory of locally compact 
groups  \cite{montgomeryzippin}; however, it still relies
on Tits' alternative for linear groups \cite{tits}
(or the easier theorem of Shalom that amenable linear
groups are virtually solvable \cite{shalom}).  

\begin{remark}
There are several important
applications of the Wilkie-Van Den Dries refinement
\cite{wilkie} of Gromov's theorem \cite{polygrowth} that do not 
follow from the original statement; for instance 
\cite{panos}, or the theorem of Varopoulos that a group
satisfies a $d$-dimensional Euclidean isoperimetric
inequality unless it is virtually nilpotent of growth 
exponent $<d$. 
\end{remark}

\subsection{Sketch of the proofs}
The proof of Theorem \ref{thmfinitedim} 
is based on a new Poincare inequality which holds for any Cayley
graph $\Ga$ of any finitely generated group $G$:  
\begin{equation}
\label{eqnmainpi}
\int_{B(R)}\;|f-f_R|^2
\leq 8\,|S|^2\,R^2\,\frac{|B(2R)|}{|B(R)|}\;
\int_{B(3R)}\;|\nabla f|^2,
\end{equation}
Here  $f$ is a piecewise smooth function on $B(3R)$,
 $f_R$ is the average of $u$ over the ball $B(R)$,
and $S$ is the generating set for $G$.

The remainder of the proof
has the same rough outline as \cite{coldingminicozzi},
though the details are different.
Note that \cite{coldingminicozzi} assumes a uniform doubling
condition as well as a uniform Poincare inequality.
In our context, we may not appeal to such uniform bounds
as their proof depends on Gromov's theorem.  Instead, the
idea is to use (\ref{eqnmainpi}) to show that one has uniform 
bounds at certain scales, and that this is sufficient
to deduce that the space of harmonic functions in 
question is finite dimensional.

\bigskip
The proof of Corollary \ref{corinfiniterep} invokes a Theorem 
of \cite{mok,korevaarschoen} to produce a fixed point
free isometric $G$-action $G\acts \h$, where $\h$ is 
a Hilbert space, and a $G$-equivariant
harmonic map $f:\Gamma\ra\h$ from the Cayley graph of $G$
to $\h$.  Theorem \ref{thmfinitedim} then implies that
$f$ takes values in a finite dimensional subspace of $\h$,
and this implies Corollary \ref{corinfiniterep}.
See Section \ref{secinfiniterep}.

\bigskip
Corollary \ref{corgromov} follows from Corollary \ref{corinfiniterep}
by induction on the degree of growth, as in the original 
proof of Gromov; see Section \ref{secgromov}.

\subsection{Acknowledgements}  I would like to
thank Alain Valette for an inspiring lecture 
at MSRI in August 2007, and the discussion afterward.
This gave me the initial impetus to find a new proof
of Gromov's theorem. 
I would especially like to thank Laurent Saloff-Coste
for telling me about the Poincare inequality in 
Theorem \ref{thmpi}, which has replaced a more complicated 
one used in an earlier draft of this paper,
and Bill Minicozzi for simplifying 
Section \ref{secproofmaintheorem}.
Finally I want to thank Toby Colding for several
conversations regarding  \cite{coldingminicozzi},
and Emmaneul Breuillard, David Fisher, Misha Kapovich, 
Bill Minicozzi, Lior Silberman and Alain Valette for 
comments and corrections.

\section{A Poincare inequality for finitely generated groups}

Let $G$ be a group, with a finite
generating set $S\subset G$.  We denote the
associated word norm of $g\in G$ by  $|g|$.
For $R\in [0,\infty)\cap \Z$, let $V(R)=|B_G(R)|=|B_G(e,R)|$.
We will denote the $R$-ball in the associated 
Cayley graph by $B(R)=B(e,R)$.  

\begin{remark} We are viewing the Cayley graph as (the
geometric realization of a) $1$-dimensional simplicial
complex, not as a discrete space.  Thus $B_G(R)$ is a finite
set, whereas $B(R)$ is typically $1$-dimensional.
\end{remark}

\begin{theorem}
\label{thmpi}
For every $R\in [0,\infty)\cap \Z$ and every smooth function
$f:B(3R)\ra \R$,
\begin{equation}
\int_{B(R)}\;|f-f_R|^2
\leq 8\,|S|^2\,R^2\,\frac{V(2R)}{V(R)}\;
\int_{B(3R)}\;|\nabla f|^2,
\end{equation}
where $f_R$ is the average of $f$ over $B(R)$.
\end{theorem}
\proof
Fix $R\in [0,\infty)\cap \Z$.

Let $\de f:B_G(3R-1)\ra \R$ be given by
$$
\de f(x)=\int_{B(x,1)}\;|\nabla f|^2.
$$

For every $y\in G$, we choose a shortest vertex
path $\ga_y:\{0,\ldots,|y|\}\ra G$ from $e\in G$ to $y$.
If $y\in B_G(2R-2)$, then
\begin{equation}
\label{eqnmult2R}
\sum_{x\in B(R-1)}\;\sum_{i=0}^{|y|}\;(\de f)(x\,\ga_y(i))
\leq 2R\;\sum_{z\in B(3R-1)}\;(\de f)(z),
\end{equation}
since the map $B(R-1)\times\{0,\ldots,|y|\}\ra B(3R-1)$
given by $(x,i)\mapsto x\,\ga_y(i)$ is at most $2R$-to-$1$.

For every ordered pair $(e_1,e_2)$ 
of edges contained in $B(R)$, let $x_i\in e_i\cap G$
be elements such that $d(x_1,x_2)\leq 2R-2$, and let
$y=x_1^{-1}x_2$. 
By the Cauchy-Schwarz inequality,
\begin{equation}
\label{eqnp1p2}
\int_{(p_1,p_2)\in e_1\times e_2}\;|f(p_1)-f(p_2)|^2\;dp_1dp_2
\;\leq\; 2R\sum_{i=0}^{|y|}\;(\de f)(x_1\,\ga_y(i)).
\end{equation}

Now
$$
\int_{B(R)}\;|f-f_R|^2
\leq \frac{1}{V(R)}\;\int_{B(R)\times B(R)}\;|f(p_1)-f(p_2)|^2\;
dp_1dp_2
$$
$$
=
\frac{1}{V(R)}\;\sum_{(e_1,e_2)\subset B(R)\times B(R)}\;
\int_{(p_1,p_2)\in e_1\times e_2}\;
|f(p_1)-f(p_2)|^2\;dp_1dp_2
$$
$$
\leq \frac{1}{V(R)}\;\sum_{(e_1,e_2)\subset B(R)\times B(R)}\;
2R\sum_{i=0}^{|y|}\;(\de f)(x_1\,\ga_y(i)),
$$
where $x_1$ and $y$ are as defined above.
The map $(e_1,e_2)\mapsto (x_1,y)$ is at most $|S|^2$-to-one,
so
$$
\int_{B(R)}\;|f-f_R|^2
\;\leq\; 2R\,|S|^2\;\frac{1}{V(R)}\;\sum_{x_1\in B(R-1)}\;
\sum_{y\in B(2R-2)}\;\sum_{i=0}^{|y|}\;(\de f)(x_1\,\ga_y(i))
$$
$$
\;\leq\; 4\,R^2\,|S|^2\frac{1}{V(R)}\;
\sum_{y\in B(2R-2)}\sum_{z\in B(3R-1)}  (\de f)(z)
\quad\mbox{(by (\ref{eqnmult2R}) )}
$$
$$
=4\,R^2\,|S|^2\,\frac{V(2R)}{V(R)}\;
\sum_{z\in B(3R-1)}  (\de f)(z)
\;\leq\; 8\,R^2\,|S|^2\,\frac{V(2R)}{V(R)}\;
\int_{B(3R)}\;|\nabla f|^2.
$$

\qed

\begin{remark}
Although the theorem above is not in the literature,
the proof is virtually contained in 
\cite[pp.308-310]{coulhonsaloffcoste}.  When hearing
of my more complicated Poincare inequality, 
Laurent Saloff-Coste's immediate response was to state
and prove Theorem \ref{thmpi}.
\end{remark}

\section{The proof of Theorem \ref{thmfinitedim}}
\label{secproofmaintheorem}

In this section $G$ will be a finitely generated group with a fixed 
 finite
generating set $S$, and the associated Cayley graph and word 
norm 
will be
denoted $\Ga$ and $\|\cdot\|$, respectively.  For $R\in \Z_+$ we 
let $B(R)\defeq B(e,R)\subset\Ga$
and $V(R)\defeq |B_G(R)|=|B(R)\cap G|$.

Let $\V$ be a $2k$-dimensional vector space of harmonic functions 
on $\Ga$.
We equip $\V$ with the family  of quadratic forms 
$\{Q_R\}_{R\in [0,\infty)}$,
where 
$$
Q_R(u,u)\defeq \int_{B(R)}\;u^2.
$$

The remainder of this section is devoted to  proving
 the following statement, which clearly implies 
Theorem \ref{thmfinitedim}:
\begin{theorem}
\label{thmfinitedim'}
For every  $d\in (0,\infty)$ there is a $C=C(d)\in (0,\infty)$ such 
that if
\begin{equation}
\label{eqnliminfbounded}
\liminf_{R\ra \infty}\;
\frac{V(R)\left(\det Q_R\right)^{\frac{1}{\dim\V}}} {R^d}\;<\;\infty,
\end{equation}
then $\dim\V<C$.
\end{theorem}

The overall structure of the proof is similar to that of 
Colding-Minicozzi \cite{coldingminicozzi}.

\subsection{Finding good scales}
We begin by using the polynomial growth assumption to select 
a pair of comparable scales $R_1<R_2$ at which both 
the growth function $V$ and the determinant 
$\left(\det Q_R\right)^{\frac{1}{\dim\V}}$ have doubling 
behavior.  Later we will use this to find many 
functions in $\V$ which have doubling behavior at scale $R_2$.  
Similar scale selection arguments appear in both 
\cite{polygrowth} and \cite{coldingminicozzi}; the one 
here is a hybrid of the two.

\bigskip
Observe that the family of quadratic
forms  $\{Q_R\}_{R\in[0,\infty)}$ is nondecreasing
in $R$, in the sense that $Q_{R'}-Q_{R}$ is positive
semi-definite when $R'\geq R$.  Also,
note that $Q_R$ is  positive definite for sufficiently 
large $R$, since $Q_R(u,u)=0$
for all $R$ only if $u\equiv 0$.  Choose $i_0\in\N$
such that $Q_R>0$ whenever $R\geq 16^{i_0}$.

We define $f:\Z_+\ra \R$ and $h:\Z\cap[i_0,\infty)\ra \R$ by
$$
f(R)=V(R)\left(\det Q_R\right)^{\frac{1}{\dim\V}},
\quad\mbox{and}\quad h(i)=\log f(16^i).$$
Note that since $Q_R$ is a nondecreasing function of $R$, 
both $f$ and $h$ are nondecreasing functions, and 
(\ref{eqnliminfbounded}) translates to
\begin{equation}
\label{eqnhbound}
\liminf_{i\ra\infty}\;\left(h(i)-di\log 16\right)<\infty.
\end{equation}

  Put $a=4d\log 16$, and pick $w\in \N$. 
\begin{lemma}
There are integers $i_1,i_2\in [i_0,\infty)$
such that 
\begin{equation}
i_2-i_1\in (w,3w),
\end{equation}
\begin{equation}
\label{eqni1i2wa}
h(i_2+1)-h(i_1)< wa,
\end{equation}
and
\begin{equation}
\label{eqni1i2a}
h(i_1+1)-h(i_1)<a, \quad h(i_2+1)-h(i_2)<a.
\end{equation}
\end{lemma}
\proof
There is a nonnegative integer $j_0$ such that 
\begin{equation}
\label{eqnj0wa}
h(i_0+3w(j_0+1))-h(i_0+3wj_0)<wa.
\end{equation}  Otherwise, for all 
$l\in \N$ we would get
$$
h(i_0+3wl)=h(i_0)+\sum_{j=0}^{l-1}\;(h(i_0+3w(j+1))-h(i_0+3wj))
$$
$$
\geq h(i_0)+wal=h(i_0)+\left(\frac43 d\log 16\right)\left(3wl\right),
$$
which contradicts (\ref{eqnhbound}) for large $l$.

Let $m\defeq i_0+3wj_0$.  

Then there are integers 
$i_1\in [m,m+w)$ and $i_2\in [m+2w,m+3w)$ such that
(\ref{eqni1i2a}) holds, for otherwise we would have
either $h(m+w)-h(m)\geq wa$ or $h(m+3w)-h(m+2w)\geq wa$,
contradicting (\ref{eqnj0wa}).  

These $i_1$ and $i_2$ satisfy the conditions of the lemma, 
because
$$
h(i_2+1)-h(i_1)\leq h(m+3w)-h(m)<wa.
$$

\qed

\subsection{A controlled cover} 
Let $R_1=2\cdot 16^{i_1}$ and $R_2=16^{i_2}$.
Choose  a maximal $R_1$-separated subset $\{x_j\}_{j\in J}$ 
of $B(R_2)\cap G$, and let $B_j\defeq B(x_j,R_1)$.  Then the
collection $\B\defeq \{B_j\}_{j\in J}$ covers $B(R_2)$,
and $\frac12\B\defeq\{\frac12 B_j\}_{j\in J}$ is a
disjoint collection.  

\begin{lemma}
\label{lembmisc}
\begin{enumerate}
\item The covers $\B$ and 
$3\B\defeq \{3B_j\}_{j\in J}$ have
 intersection multiplicity $<e^a$.

\item $\B$ has cardinality $|J|<e^{wa}$.

\item There is a $C\in (0,\infty)$ depending only
on $|S|$ such that
for every  $j\in J$ and every smooth function 
$v:3B_i\ra \R$,
\begin{equation}
\int_{B_i}\;|v-v_{B_i}|^2\leq C\,e^a\,R_1^2
\int_{3B_i}\;|\nabla v|^2.
\end{equation} 
\end{enumerate}
\end{lemma}
\proof

(1)  If $z\in 3B_{j_1}\cap\ldots\cap 3B_{j_l}$,
then $x_{j_m}\in B(x_{j_1},6R_1)$ for every $m\in \{1,\ldots,l\}$,
so $\{B(x_{j_m},\frac{R_1}{2})\}_{m=1}^l$ are disjoint balls
lying in $B(x_{j_1},8R_1)$, and hence 
$$
\log l\leq \log\frac{V(3R_1)}{V(\frac{R_1}{2})}
=\log V(3R_1)-\log V\left(\frac{R_1}{2}\right)
\leq h(i_1+1)-h(i_1)<a.$$
This shows that the multiplicity of $3\B$ is  at most
$e^a$.  This implies (1), since the multiplicity of 
$\B$ is not greater than that of $3\B$.

\bigskip
(2)  The balls $\{B(x_j,\frac{R_1}{2})\}_{j\in J}$
are disjoint, and are contained in $B(R_2+\frac{R_1}{2})\subset
B(2R_2)$, so 
$$
|J|\leq \frac{V(2R_2)}{ V(\frac{R_1}{2})}\leq 
\frac{V(16^{i_2+1})}{V(16^{i_1})}< e^{wa},$$
by (\ref{eqni1i2wa}).

\bigskip
(3)  By Theorem \ref{thmpi} and the translation invariance
of the inequality, 
$$
\int_{B_i}\;|v-v_{B_i}|^2
\;\leq\; 8\,|S|^2\,R_1^2\frac{V(2R_1)}{V(R_1)}
\int_{3B_i}\;|\nabla v|^2
$$
$$
\leq 8\,|S|^2\,R_1^2\,e^a
\int_{3B_i}\;|\nabla v|^2.
$$
\qed

\subsection{Estimating functions relative to the cover $\B$}
We now estimate the size of a  harmonic function
in terms of its averages over the $B_j$'s, and its 
size on a larger ball.

We define a linear map $\Phi:\V\ra \R^J$ by
$$
\Phi_j(v)\defeq \frac{1}{|B_j|}\int_{B_j}\;v.
$$

\begin{lemma}[cf.\protect{ \cite[Prop. 2.5]{coldingminicozzi}}]
There is a constant $C\in(0,\infty)$ depending only
on the size of the generating set $S$, with the following property. 

\begin{enumerate}

\item
If $u$ is a smooth functions on $B(16R_2)$, then
\begin{equation}
\label{eqnnonharmonicbound}
Q_{R_2}(u,u)
\quad\leq\quad C\,V(R_1)\,|\Phi(u)|^2
+\;C\,e^{2a}R_1^2\;
\int_{B(2R_2)}\;|\nabla u|^2.
\end{equation}

\item  If $u$ is harmonic on $B(16R_2)$, then
\begin{equation}
\label{eqnharmonicbound}
\begin{aligned}
Q_{R_2}(u,u)
\leq C\,V(R_1)|\Phi(u)|^2+
C\,e^{2a}\left(\frac{R_1}{R_2}\right)^2
Q_{16R_2}(u,u).
\end{aligned}
\end{equation}
\end{enumerate}
\end{lemma}
\proof
We will use $C$ to denote a constant which depends 
only on $|S|$; however, its value may vary from equation to 
equation.

We have
$$
Q_{R_2}(u,u)=\int_{B(R_2)}\;u^2
\leq \sum_{j\in J}\;\int_{B_j}\;u^2
$$
\begin{equation}
\label{eqn3term}
\leq 2\sum_{j\in J}\;\int_{B_j}\;
\left(|\Phi_j(u)|^2+
|u-\Phi_j(u)|^2\right).
\end{equation}

We estimate each of the terms in (\ref{eqn3term})
in turn.

For the first term we get:
\begin{equation}
\label{eqnfirstterm}
\sum_{j\in J}\;\int_{B_j}\;|\Phi_j(u)|^2
=\sum_{j\in J}\;|B_j|\left|\Phi_j(u)\right|^2
\leq C\,V(R_1)\;|\Phi(u)|^2.
\end{equation}

For the second term we have:
$$
\sum_{j\in J}\;\int_{B_j}\;|u-\Phi_j(u)|^2
$$
$$
\leq C\,e^a\,R_1^2\,\sum_{j\in J}\;
\int_{3B_j}\;|\nabla u|^2
\quad\mbox{(by (3) of Lemma \ref{lembmisc})}
$$
$$
\leq Ce^aR_1^2\;\left(e^a
\int_{B(2R_2)}\;|\nabla u|^2\right)
\quad\mbox{(by (1) of Lemma \ref{lembmisc})}
$$
$$
=Ce^{2a}\,R_1^2\;
\int_{B(2R_2)}\;|\nabla u|^2.
$$
 Combining this with 
(\ref{eqnfirstterm}) yields (1).

\bigskip
Inequality (\ref{eqnharmonicbound}) follows from 
(\ref{eqnnonharmonicbound}) by applying the
reverse Poincare inequality, which holds for any
harmonic function $v$ defined on $B(16R_2)$:
$$
R_2^2\int_{B(2R_2)}\;|\nabla v|^2
\leq C\;Q_{16R_2}(v,v).
$$
(For the proof, see \cite[Lemma 6.3]{schoenyau}, and note 
that for harmonic
functions their condition $u\geq 0$ may be dropped.)
\qed

\bigskip
\bigskip

\subsection{Selecting functions from $\V$ with controlled
growth}
Our next step is to select functions in $\V$ which have doubling
behavior at scale $R_2$.  

\begin{lemma}[cf.\protect{ \cite[Prop. 4.16]{coldingminicozzi}}]
\label{lemselectu}
There is a subspace $\U\subset\V$ of dimension at
least $k=\frac{\dim\V}{2}$ such that for every $u\in \U$
\begin{equation}
\label{eqnR2controls16R2}
Q_{16R_2}(u,u)\leq e^{2a}\;Q_{R_2}(u,u).
\end{equation}

\end{lemma}
\proof
Since $R_2=16^{i_2}>16^{i_0}$, the quadratic form $Q_{R_2}$
is positive definite.  Therefore there is a $Q_{R_2}$-orthonormal
basis $\be=\{v_1,\ldots,v_{2k}\}$ for $\V$ which is orthogonal
with respect to $Q_{16R_2}$.   

Suppose there are at least $l$ distinct elements $v\in\be$ such that 
 $Q_{16R_2}(v,v)\geq e^{2a}$.  Then since $\be$ is
$Q_{R_2}$-orthonormal and $Q_{16R_2}$-orthogonal,
$$
\log\left(\frac{\det Q_{16R_2}}{\det Q_{R_2}}\right)^\frac{1}{2k}
=\log\left(\prod_{j=1}^{2k}\frac{Q_{16R_2}(v_j,v_j)}{Q_{R_2}(v_j,v_j)}
\right)^{\frac{1}{2k}}
=\log\left(\prod_{j=1}^{2k}\;Q_{16R_i}(v_j,v_j)\right)^{\frac{1}{2k}}
$$
$$
\geq \log\left(e^{2al}\right)^{\frac{1}{2k}}=\frac{l}{k}a.
$$
On the other hand,
$$
a> h(i_2+1)-h(i_2)\geq \log\left(\det Q_{16R_2}\right)^{\frac{1}{2k}}
-\log\left(\det Q_{R_2}\right)^{\frac{1}{2k}}.
$$
So we have a contradiction if $l\geq k$.  

Therefore we may
choose a $k$ element subset 
$\{u_1,\ldots,u_k\}\subset\{v_1,\ldots,v_{2k}\}$ such that
$Q_{16R_2}(u_j,u_j)<e^{2a}$ for every $j\in \{1,\ldots,k\}$.
Then every element of $\U\defeq\span\{u_1,\ldots,u_k\}$ 
satisfies (\ref{eqnR2controls16R2}). 
\qed

\mbox{}
\bigskip
\bigskip

\subsection{Bounding the dimension of $\V$}
We now assume that 
$w$ is the smallest integer such that
\begin{equation}
\label{eqnwbig}
\left(\frac{R_1}{R_2}\right)^2=2\cdot 16^{i_1-i_2}<2\cdot 16^{-w}<
\frac{1}{2Ce^{4a}},
\end{equation}
where $C$ is the constant in (\ref{eqnharmonicbound}).  Therefore
$2\cdot 16^{-(w-1)}\geq \frac{1}{2Ce^{4a}}$, and this implies
\begin{equation}
\label{eqnewaestimate}
e^{wa}\leq 64\,C\,e^{64d^2\log 16}.
\end{equation}

If $u\in \U$ lies in the kernel of $\Phi$, then 
$$
Q_{R_2}(u,u)\leq
Ce^{2a}\left(\frac{R_1}{R_2}\right)^2\,
Q_{16R_2}(u,u)\quad\mbox{(by (\ref{eqnharmonicbound}) )}
$$
$$
\leq 
Ce^{2a}\left(\frac{R_1}{R_2}\right)^2\,\left(e^{2a}\;
Q_{R_2}(u,u)\right)\quad\mbox{(by  Lemma \ref{lemselectu})}
$$
$$
\leq\frac12 Q_{R_2}(u,u)\quad\mbox{(by (\ref{eqnwbig}) )}.
$$
Therefore $u=0$, and we conclude that
$\Phi\restr_{\U}$ is injective.  Hence  by Lemma
\ref{lembmisc}
and (\ref{eqnewaestimate}),
$$
\dim\V=2\dim\U\leq 2|J|\leq 2e^{wa}\leq 128\,C\,e^{64d^2\log 16}.
$$

\qed

\section{Proof of Corollary \ref{corinfiniterep} using 
Theorem \ref{thmfinitedim}}
\label{secinfiniterep}

Let $G$ be as in the statement of the Corollary, 
and let $\Ga$ denote some Cayley graph of $G$ with 
respect to a symmetric finite generating set $S$. 

Note that $G$ is amenable, for if $R_k\ra \infty$
and $V(R_k)<AR_k^d$ for all $k$, then for every $k$
there must be an $r_k\in [\frac{R_k}{2},R_k]$ such
that the ball $B_G(r_k)$ satisfies
$$
|\D B_G(r_k)|=|S_G(r_k)|<3A\,R_k^{d-1};
$$
this means that the sequence of 
balls  $\{B_G(r_k)\}$ is a Folner sequence for $G$.  

Hence $G$  does not have Property (T).  
Therefore
by a result of Mok \cite{mok}
and Korevaar-Schoen \cite[Theorem 4.1.2]{korevaarschoen}, 
there is an isometric 
action $G\acts \h$ of $G$ on a 
Hilbert space $\h$ which has no fixed points,  and a nonconstant 
$G$-equivariant 
harmonic map $f:\Ga\ra \h$.  In the case of
Cayley graphs, the Mok/Korevaar-Schoen result is quite elementary,
so we give a short proof in Appendix A.

Since $f$ is $G$-equivariant,
it is Lipschitz.

Each bounded linear functional $\phi\in\h^*$ gives rise to a 
Lipschitz harmonic
function $\phi\circ f$, and hence we have a linear map 
$\Phi:\h^*\ra \V$, where $\V$ is the space of Lipschitz
harmonic functions on $\Ga$.  Since the target is finite dimensional
by Theorem \ref{thmfinitedim}, the kernel of $\Phi$ has finite
codimension, and its annihilator $\ker(\Phi)^\perp\subset\h$ is
a  finite dimensional subspace containing the image of $f$.
It follows that the affine hull $A$ of the image of $f$ is 
finite dimensional
and $G$-invariant.  Therefore we have an induced isometric 
$G$-action $G\acts A$.  
This action cannot factor through a finite group, because it 
would then have 
fixed points, contradicting the fact that the original 
representation
is fixed point free.  The associated homomorphism 
$G\ra \isom(A)$ yields the
desired finite dimensional representation of $G$.

\qed

\section{Proof of Corollary \ref{corgromov} using Corollary 
\ref{corinfiniterep}}
\label{secgromov}

We prove Gromov's theorem using Corollary \ref{corinfiniterep}.
The proof is a recapitulation of Gromov's argument, which  
reproduce here for the convenience of the reader.

The proof is by induction on the degree of growth.

\begin{definition}
Let $G$ be a finitely generated group.  The {\bf degree (of growth)
of $G$} is the minimum $\deg(G)$ of the  nonnegative
integers $d$ such that 
$$
\liminf_{r\ra\infty}\;\frac{V(r)}{r^d}<\infty.
$$
\end{definition}
A group whose degree of growth is $0$ is finite, and hence
Corollary \ref{corgromov} holds for such a group.

Assume inductively that for some $d\in \N$
that every group of degree at most $d-1$ is virtually
nilpotent, and suppose $\deg(G)=d$.   Then $G$ is infinite,  
and by Corollary \ref{corinfiniterep}  there is a finite 
dimensional
linear representation $G\ra GL(n)$ with infinite image $H\subset GL(n)$.
Since $H$ has polynomial growth, by \cite{tits} (see \cite{shalom}
for an easier proof) it is virtually solvable,
and by \cite{wolf,milnor} it must be virtually nilpotent.  

After passing to 
finite index subgroups, we may assume $H$ is nilpotent, and that its 
abelianization
is torsion-free.  It follows that there is a short exact sequence
$$
1\lra K\ra G\stackrel{\al}{\ra} \Z\lra 1.
$$
By \cite[Lemma (2.1)]{wilkie}, the normal subgroup $K$ is 
finitely generated,
and $\deg(K)\leq \deg(G)-1$.  

By the induction hypothesis, $K$ is virtually nilpotent. 
Let $K'$ be a 
finite index nilpotent subgroup of $K$ which is normal in $G$, 
and let
$L\subset G$ be an infinite cyclic subgroup which
is mapped  isomorphically by $\al$ onto $\Z$.  
Then $K'L\subset G$ 
is a finite
index solvable subgroup of $G$.  As it has polynomial growth, by 
\cite{wolf,milnor} it is virtually nilpotent.

\qed

\appendix

\section{Property (T) and equivariant harmonic maps}

In this expository section, we will give a simple proof of the 
special case of the Korevaar-Schoen/Mok existence result needed 
in the proof of Corollary \ref{corgromov}.

Suppose $G$ is a finitely generated group, $S=S^{-1}\subset G$
is a symmetric finite generating set, and $\Ga$
is the associated Cayley graph.

Given an  action $G\acts X$ on a metric space $X$, we
define the {\bf energy function} $E:X\ra \R$ by 
$$
E(x)=\sum_{s\in S}\;d^2(s x,x).
$$

We recall that a $G$ has Property (T) iff every isometric
action of $G$ on a Hilbert space has a fixed point.

The following theorem is a very weak version
of some results in  \cite{fishermargulis}, see also
\cite[pp.115-116]{randomgroups}:

\begin{theorem}
\label{thma1}
The following are equivalent:
\begin{enumerate}
\item $G$ has Property (T).
\item There is a constant $D\in (0,\infty)$ such that
if $G\acts\h$ is an isometric action on a Hilbert space
and $x\in \h$, then $G$ fixes a point in $B(x,D\sqrt{E(x)})$.
\item There are constants $D\in (0,\infty)$, $\la\in (0,1)$
such that if $G\acts \h$ is an isometric action on 
a Hilbert space and $x\in \h$, then there is a point
$x'\in B(x,D\sqrt{E(x)})$ such that $E(x')\leq \la E(x)$.
\item There is no isometric action $G\acts \h$ on a Hilbert
space such that the energy function $E:\h\ra\R$ attains a
positive minimum.
\end{enumerate}
\end{theorem}
\proof
Clearly (2)$\implies$(1).  Also,  (1)$\implies$(4) since
the energy function $E$ is zero at a fixed point.

\bigskip
(3)$\implies$(2).  Suppose (3) holds.  Let $G\acts \h$
be an isometric action, and pick $x_0\in \h$.  Define a 
sequence $\{x_k\}\subset \h$ inductively, by choosing
$x_{k+1}\in B(x_k,D\sqrt{E(x_k)})$ such that 
$E(x_{k+1})\leq \la E(x_k)$.  Then $E(x_k)\leq \la^k\,E(x_0)$
and $d(x_{k+1},x_k)\leq D\sqrt{E(x_k)}
\leq D\la^{\frac{k}{2}}\sqrt{E(x_0)}$.  Therefore
$\{x_k\}$ is Cauchy, with limit $x_\infty$ satisfying
 $$
d(x_\infty,x_0)\leq \frac{D\sqrt{E(x_0)}}{1-\la^{\frac12}}.
$$
Then $E(x_\infty)=\lim_{k\ra\infty}\;E(x_k)=0$, and
$x_\infty$ is fixed by $G$.  Therefore (2) holds.

\bigskip
(4)$\implies$(3).  We prove the contrapositive.  
Assume that (3) fails.  Then for every $k\in \N$,
we can find an isometric action $G\acts \h_k$ on 
a Hilbert space, and a point $x_k\in \h_k$ such that
\begin{equation}
\label{eqnenergy}
E(y)>\left(1-\frac{1}{k}\right)\,E(x_k)
\end{equation}
for every $y\in B(x_k,k\sqrt{E(x_k)})$.  Note that in 
particular, $E(x_k)>\left(1-\frac{1}{k}\right)\,E(x_k)$,
forcing $E(x_k)>0$.

Let $\h_k'$ be the result of rescaling the metric on $\h_k$
by $\frac{1}{\sqrt{E(x_k)}}$.   Then (\ref{eqnenergy})
implies that the induced isometric
action $G\acts\h_k'$ satisfies $E(x_k)=1$ and 
\begin{equation}
E(y)\geq 1-\frac{1}{k}
\end{equation}
for all $y\in B(x_k,k)$.  Then any ultralimit (see 
\cite{asyinv,kleinerleeb})
of the 
sequence $(\h_k,x_k)$ of pointed Hilbert spaces is
a pointed Hilbert space $(\h_\om,x_\om)$ with an 
isometric action $G\acts \h_\om$ such that 
$$
E(x_\om)=1=\inf_{y\in \h_\om}\;E(y).
$$
Therefore (4) fails.

\qed

\bigskip
\bigskip
Before proceeding we recall some facts about harmonic
maps on graphs.   Suppose $\G$ is a locally finite metric
graph, where all edges have length $1$.   If
$f:\G\ra \h$ is a piecewise
smooth map to a Hilbert space, then the following
are equivalent:

\begin{itemize}
\item $f$ is harmonic.
\item The Dirichlet energy of $f$ (on any finite subgraph)
is stationary with respect to compactly supported variations
of $f$.
\item The restriction of $f$ to each edge of $\G$
has constant derivative, and for every vertex $v\in \G$, 
$$
\sum_{d(w,v)=1}\;\left(f(w)-f(v)\right)=0.
$$ 
\end{itemize}

Note that if $G\acts \h$ is an isometric action on a Hilbert
space, then $E$ is a smooth convex function,  and its derivative
is
$$
DE(x)(v)=2\left(\sum_{s\in S}\;\langle s x-x,(Ds)( v)\rangle
-\sum_{s\in S}\;\langle s x-x,v\rangle\right)
$$
$$
=2\left(\sum_{s\in S}\;\langle x-s^{-1}x,v\rangle
+\sum_{s\in S}\;\langle x-s x,v\rangle\right)
=4\sum_{s\in S}\;\langle x-s x,v\rangle.
$$
Therefore  
$$
\mbox{$x\in\h$ is a critical point of $E$}
$$
$$
\Longleftrightarrow\mbox{ $x$
is a minimum of $E$}
$$
\begin{equation}
\label{eqncriticalpoint}
\Longleftrightarrow
\sum_{s\in S}\;(x-s x)=0.
\end{equation}
Therefore the $G$-equivariant map $f_0:G\ra \h$
given by $f_0(g)\defeq gx$ extends to a $G$-equivariant
harmonic  map $f:\Ga\ra\h$ if and only if
$$
\sum_{s\in S}\;\left(f_0(s e))-f_0(e)\right)
=\sum_{s\in S}\;\left(s x-x\right)=0
$$
$$
\Longleftrightarrow\mbox{$x$ is a minimum of $E$.}
$$

\bigskip
The next result is a very special case of a theorem
from \cite{mok,korevaarschoen}.
\begin{lemma}
The following are equivalent:
\begin{enumerate}
\item
$G$ does not have Property (T).
\item There is an isometric
action $G\acts \h$ on a Hilbert space $\h$ and a
nonconstant $G$-equivariant harmonic map $f:\Ga\ra\h$.
\end{enumerate}
\end{lemma}
\proof
(1)$\implies$(2).  If $G$ does not have Property (T), then
by Theorem \ref{thma1} there is an isometric action 
$G\acts \h$ on a Hilbert space,  and a point $x\in \h$ with 
$E(x)=\inf_{y\in \h}\;E(y)>0$.  Let $f:\Ga\ra \h$ be the
 $G$-equivariant map with $f(g)=gx$ for every $g\in G\subset\Ga$, 
and whose
restriction to each edge $e$ of $\Ga$ has constant derivative.
Then $f$ is harmonic,
and obviously nonconstant.

\bigskip
(2)$\implies$(1).  Suppose (2) holds, and $f:\Ga\ra\h$
is the $G$-equivariant harmonic map.  Then 
 $f(e)$ is a positive minimum of $E:\h\ra\R$;
in particular the action $G\acts\h$
has no fixed points.  Therefore $G$ does not have 
Property (T). 

\qed

\bibliography{polygrowth}
\bibliographystyle{alpha}

\end{document}